\documentclass[pdftex,12pt,a4paper]{article}
\pdfoutput=1
\usepackage{graphicx}
\usepackage{float}
\usepackage{wrapfig}
\usepackage[rightcaption]{sidecap}
\usepackage{subfig}
\usepackage{url}
\usepackage{hyperref}
\usepackage{placeins}
\usepackage{amsmath}

\graphicspath{ {images/} }
\DeclareMathSizes{14}{24}{12}{10}
\usepackage[margin=1in]{geometry}

\begin{document}

{\centering
{\textbf{\Large An Implementation of Adaptive Mesh Refinement for Shallow Water Equations}}
\medskip\\
\centerline{\large Avi Schwarzschild\footnote{Department of Applied Mathematics, University of Washington, Seattle, WA. (avi1@uw.edu)} and Kyle Mandli\footnote{Department of Applied Physics and Applied Mathematics, Columbia University, New York, NY.} }}

\medskip

\begin{abstract}
An implementation of adaptive mesh refinement algorithms is presented for use with multilayer shallow water equations. Currently, adaptive mesh refinement is implemented with a single layer shallow water model in the GeoClaw framework. This implementation, also in the GeoClaw framework, is for multilayer models, which have been implemented in GeoClaw previously. Until now, however, these models were too computationally expensive to run on large domains while resolving detail in coastal regions.
\end{abstract}

\section{Multilayer Adaptive Mesh Refinement}

Tsunami and storm surge are often modeled with the two-dimensional depth-averaged shallow water equations (SWE) equations. While the SWE use a single layer model, a multilayer approach may better capture certain phenomena, for example wind driven storm surge \cite{thesis}. The GeoClaw software \cite{BergerGeorgeLeVequeMandli:awr11,an,mandli2014adaptive}, a part of Clawpack designed to solve problems involving shallow water flows, includes solvers for the multilayer SWE \cite{mandli}. One of the major difficulties in the numerical simulation of any wave propagation problem on a large domain is computing results that capture the detail in the solution in a reasonable amount of time. Adaptive mesh refinement (AMR) is a method that has been used to manage the resolution of the discretization of the domain to provide detail where it is important, without increasing the number of computations necessary. Geoclaw has previously had AMR implemented with single layer SWE models \cite{an}. Using the multilayer model on geophysical problems without AMR is unfeasible at a large scale. We describe new refinement and coarsening algorithms which have been implemented for use with the multilayer shallow water equations in GeoClaw, and enable solutions to be computed. We also present an example calculation with timing data, to demonstrate the difference in speed when these algorithms are used.

\section{AMR}
A comprehensive overview of the AMR algorithm used in GeoClaw can be found in \cite{an} and, given the similarity, only a brief overview will be given here. The general idea is that the grid on which the solution is computed can adaptively vary in resolution as the simulation progresses. Small $\Delta x$ and $\Delta y$ are necessary in  regions where the solution is varying, i.e. near larger waves and inundation areas. These dense grids require many computations by the numerical solver. Therefore, it is desirable to have coarser grids in much of the domain, where good results can be achieved even with lower resolution. The management of these finer and coarser meshes is adaptive as it follows details in the solution. As a result, at a given location, we need algorithms to move from one mesh level to another. The two directions, refinement and coarsening, correspond to different sets of algorithms, both are further explained below.\\

For a given mesh level $\ell$ the spacial refinement ratios $r^\ell_x$ and $r^\ell_y$ dictate how many cells fit in one cell of the next coarser level. 
\begin{equation}
    \Delta x^{\ell+1} = \Delta x^\ell / r^\ell_x, 
\end{equation}
\begin{equation}
    \Delta y^{\ell+1} = \Delta y^\ell / r^\ell_y.
\end{equation}
Often, in practice $r_x^\ell = r_y^\ell$. Additionally, in order to maintain a Courant number less than one, there must be a temporal refinement factor $r_t^\ell$. When solving hyperbolic PDEs with explicit methods, in general $r_t^\ell = \max(r_x^\ell, r_y^\ell)$ to respect the CFL condition. However, if fine grids are used only in shallow regions, where the shallow water wave speed is smaller, it may be possible to use smaller $r_t^\ell$ and GeoClaw chooses this ratio dynamically \cite{an}.

\subsection{Refinement}

The AMR procedures to create new mesh grids of higher resolution are similar to the single layer algorithms with two major changes. The first is in the interpolation algorithm, the second is an added iterative loop to march through all the layers. We begin this section by describing the modified interpolation scheme. For simplicity, we will show eqauations in one dimenssion. Let $m \in \{1,2, ... M\}$ index the layers and $k \in \{1,2,3...\}$ index the cells on mesh level $\ell$. The layers are such that layer $m=1$ is the top layer, and any layer can only be wet in a given cell if all the layers above it in that location are also wet. 

Initializing new grids with finer resolution requires refining the bathymetry as well as interpolating the surface elevations $\eta^\ell_m$ of each layer. We refine the surface elevations rather than the fluid depths $h_m^\ell$ in order to conserve mass and prevent the generation of new surface extrema. It is important to note that the refinement algorithm presented here and the one in \cite{an} will not perfectly conserve mass when refining near the shore. This cannot be avoided when the topography is resolved on different resolution grids. A longer discussion on the topic can be found in \cite{an,thesis,mandli}.

Previously within a given cell: 
\begin{equation}
 \eta = h + B,    
 \label{eq:hb}
\end{equation}
where $h$ is the fluid depth, and B is the bathymetry function. In the multilayer case we need to consider the depths of all the layers beneath the surface we are interpolating. In the multilayer case the surface elevation $\eta^\ell_{m,k}$ of a layer $m$, at some mesh level $\ell$, in some cell $k$ can be computed as:
\begin{equation}
    \eta^\ell_{m,k} = B^\ell_k + h^\ell_{M,k} + h^\ell_{M-1,k} + ... + h^\ell_{m,k}
\end{equation}
We can introduce a new variable here called the effective bathymetry $\hat{B}$:
\begin{equation}
    \hat{B}_{m,k}^\ell = \sum_{n = m+1}^M h_{n,k}^\ell + B_k^\ell
\end{equation}
This will allow us to have an expression in form of \eqref{eq:hb} for each layer. For ease of describing the algorithm we will restrict ourselves to two layers, and the mess of indices above simplifies. For two layers we have two surfaces,
\begin{equation}
    \eta^\ell_{1,k} = B^\ell_k + h^\ell_{2,k} + h^\ell_{1,k}
\end{equation}
\begin{equation}
    \eta^\ell_{2,k} = B^\ell_k + h^\ell_{2,k}, 
\end{equation}
which can be rewritten using $\hat{B}$ as
\begin{equation}
    \eta^\ell_{1,k} = \hat{B}^\ell_k + h^\ell_{1,k}
\end{equation}
\begin{equation}
    \eta^\ell_{2,k} = B^\ell_k + h^\ell_{2,k}.
\end{equation}
With only two layers, the layer index on $\hat{B}$ can be ignored, since the second layer has $\hat{B}_{2,k}^\ell = B_k^\ell$. To compute $\eta$ values for the next finer mesh level we use linear interpolants. The slope $\sigma_k^\ell$, is computed as follows:

\begin{equation}
    \sigma_k^\ell = \textrm{minmod}\left(\eta_k^\ell - \eta_{k-1}^\ell,\eta_{k+1}^\ell - \eta_{k}^\ell \right) / \Delta x^\ell
\end{equation}

The minmod function here returns the minimum modulus of the two arguments, unless they have opposite signs, in which case it returns zero. This prevents the generation of new extrema and is discussed further in \cite{leveque} \cite{an} . 

The second major change is fairly intuitive, we now need to march through the layers so we have an added iterative loop to do so. This algorithm starts with interpolating the surface of the bottom layer (sitting immediately above the bathymetry), exactly like the single layer case, and then continuing to move up the layers using the effective bathymetry defined above.

\subsection{Coarsening}
When grids of more than one resolution are present, the solutions on the finer grids are used to update the coarse solutions. There is almost no difference between the way this is done in multilayer simulations and single layer simulations. In both cases, the average of the finer cells is used to update the coarser ones. Let $i \in \Gamma_k^\ell$ index the fine cells on mesh level $\ell+1$ within coarse cell $k$. Then the fluid depths on mesh level $\ell$ are computed as
\begin{equation}
h^\ell_{m,k} = \frac{1}{r^\ell_x}\sum_{i \in \Gamma^\ell_k}h^{\ell+1}_{m,i}
\end{equation}
All the quantities in the solution vector are averaged in the same way and cast down to the coarser grids.

\section{Demonstration}
This example is an idealized rectangular domain with a plane wave moving from deep water onto a shelf modeled as a discontinuous jump in the bathymetry. The dimensions here are scaled to show sea level at 0 and the deep ocean floor at -1. The images below are plots of the surface elevation with the AMR grid levels overlaid.

We tested this example with and without AMR, and shared memory parallelism using OpenMP, which is provided in GeoClaw. The timing data from these experiments are in the table below and are averaged over several runs with identical setups. The setups can be found on GitHub.\footnote{\url{https://github.com/aks2203/geoclaw/tree/timetests/examples/multi-layer/tests_to_run}} All the times are reported here in seconds. There are several important conclusions to be drawn from the table. The first is that, as expected, there is a small overhead cost associated with running these simulations in parallel. There is also overhead associated with the AMR code iteself. In the example presented here, using AMR reduced the number of cell updates to 7.20\% of the number of without AMR. The runtime was reduced to 9.57\% of the runtime without AMR. This overhead is associated with executing the refinement code, managing different resolution grids, and integrating on different levels. There is no general approach by which the speed up due to AMR can be predicted, as it is largely problem dependent. With large domains and small areas of interest AMR will help tremendously. In a small domain where the waves being tracked fill the entire domain, AMR will have a smaller impact.\\

\begin{figure}[]
\vspace*{-2cm}\hspace*{-3cm}\subfloat[]{\includegraphics[width = 1.4\textwidth]{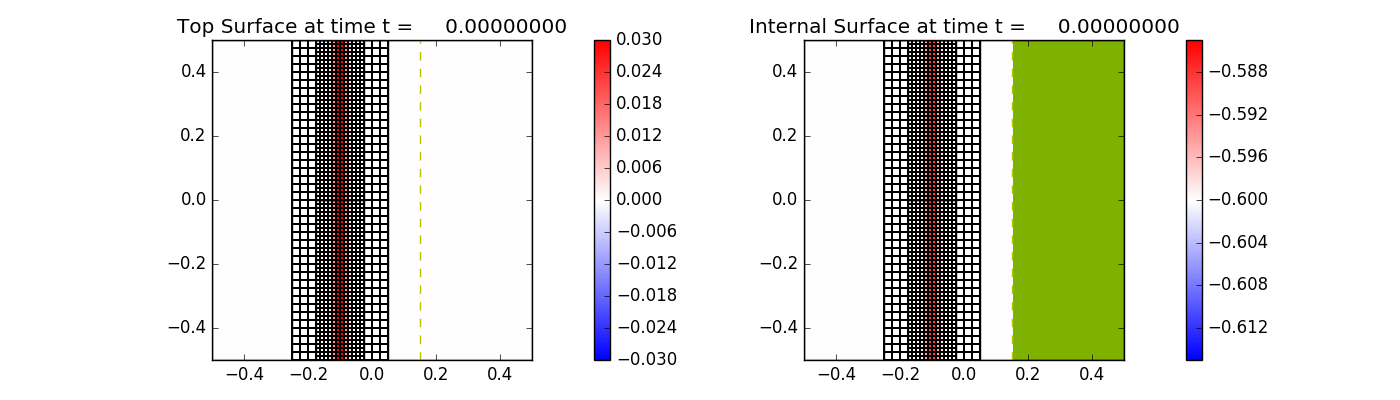}} \\
\hspace*{-3cm}\subfloat[]{\includegraphics[width = 1.4\textwidth]{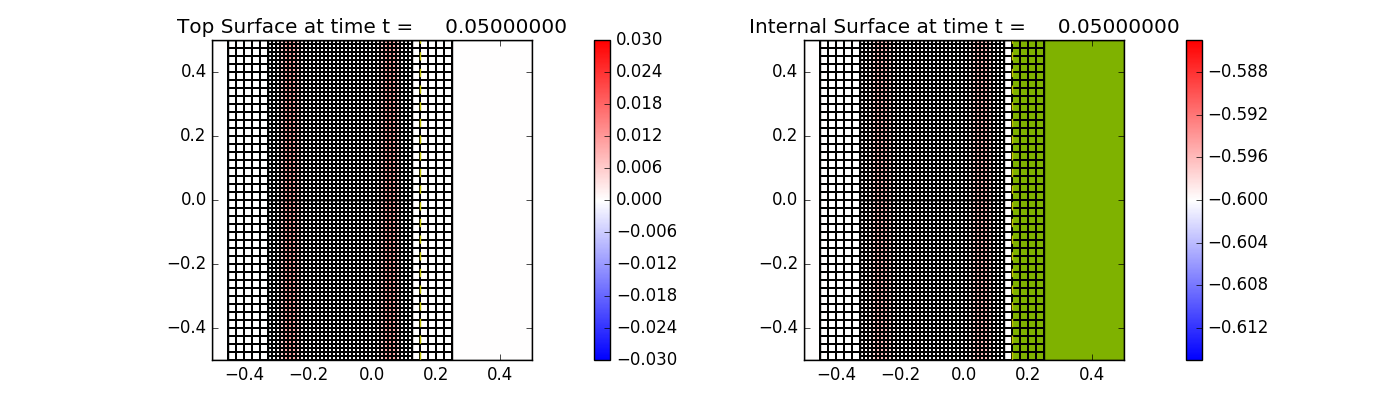}} \\
\hspace*{-3cm}\subfloat[]{\includegraphics[width = 1.4\textwidth]{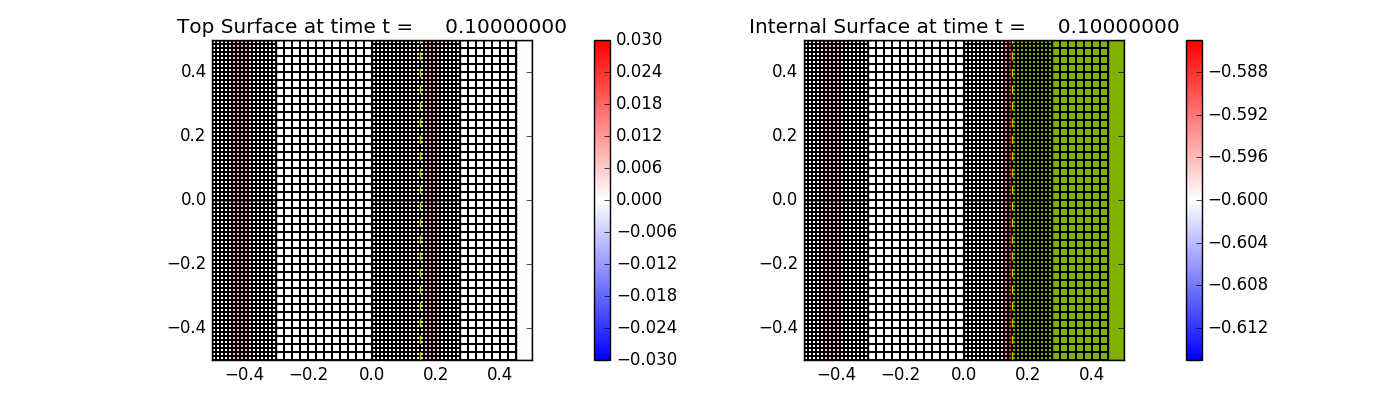}} 
\caption{Surface elevation at three different times. The images on the left show the elevation of the free surface of the ocean. The images on the right are the elevation of the internal surface between the two layers. In each plot, the grid lines show different regions of refinement. Note that the refinement follows the waves as they move accross the domain.}
\end{figure}

\begin{table}[h]
\centering
\caption{Timing data (all times in seconds)}
\begin{tabular}{c|c|c|c|c}

AMR & Threads & Avg. Wall Time & Avg. CPU Time & Total Cell Updates\\
\hline
\hline
No  &  1      &  17,580                       &    17,580    &   $2.46\times10^9$   \\        
Yes &  1      &  1,680                        &    1,680     &   $1.77\times10^8$ \\  
\hline
\hline
No  &  4      &  5,640                        &    22,500    & $2.46\times10^9$  \\   
Yes &  4      &  540                         &    2,160      & $1.77\times10^8$ \\ 

\end{tabular}
\end{table}

\section{Conclusion}
The examples run in this study confirmed that the AMR algorithms implemented in GeoClaw's multilayer models shorten the run times considerably. We showed that for 7.20\% of the number of updates, adaptively refined runs require $9.57\%$ of the amount of time. This demonstrates that the overhead cost is far outweighed by the benefit. In the future, we hope to use the multilayer models with AMR on real topography and with initial conditions to match meteorological and geological events. There is more to be explored in the parallelization of these algorithms and implementing this code for use with GPUs. 

\pagebreak
\bibliography{report} 
\bibliographystyle{abbrv}

\end{document}